\documentclass{amsart}

\newcommand{\nc}{\newcommand}

\newcommand{\fa}{{\frak a}}
\newcommand{\ff}{{\frak f}}
\newcommand{\fg}{{\frak g}}
\newcommand{\fh}{{\frak h}}

\newcommand{\fl}{{\frak l}}
\newcommand{\fn}{{\frak n}}
\nc{\fo}{{\frak o}}

\newcommand{\fs}{{\frak s}}
\newcommand{\fu}{{\frak u}}
\newcommand{\fv}{{\frak v}}

\newcommand{\fz}{{\frak z}}

\newcommand{\cc}{{\Bbb C}}

\nc{\nn}{{\Bbb N}}
\newcommand{\rr}{{\Bbb R}}
\nc{\zz}{{\Bbb Z}}

\newtheorem{teo}{Theorem}[section]

\newtheorem{prop}[teo]{Proposition}

\nc{\lp}{\langle}
\nc{\rp}{\rangle}
\nc{\inc}{\hookrightarrow}

\title[Abelian complex structures]{\bf Abelian complex structures on 
solvable Lie algebras}
\author{Mar\'{\i}a L. Barberis}
\date{}
\address{CIEM, FaMAF, Universidad Nacional de C\'ordoba, Ciudad Universitaria, (5000) C\'ordoba, Argentina}
\email{barberis@mate.uncor.edu}
\thanks{Both authors were partially supported 
by CONICET and Secyt-UNC (Argentina).}

\author{Isabel G. Dotti}
\email{idotti@mate.uncor.edu}

\subjclass{Primary 17B30, 32M10; Secondary  53C15}
\keywords{abelian complex structures, solvable Lie algebras, affine Lie algebras}

\begin{document}
\maketitle

\begin{abstract} We 
obtain a characterization of the Lie algebras admitting abelian 
complex structures in terms of 
 certain affine Lie algebras $\fa \ff \ff (A)$, where 
$A$ is a commutative algebra. 
\end{abstract}

\section{Introduction}

An {\it abelian} complex structure on a real Lie algebra $\fg$ is an endomorphism 
of $\fg$ satisfying \begin{equation} J^2=-I, 
\hspace{1.5cm} [Jx,Jy]=[x,y], \; \;\; \forall x,y 
\in \fg.   \label{abel}  \end{equation}
If $G$ is a Lie group with Lie algebra $\fg$ these conditions  imply the vanishing 
of the Nijenhuis tensor on the invariant almost complex manifold $(G,J)$, that is, $J$ is integrable on $G$. 

Our interest arises from properties of the complex manifolds obtained by considering this class 
of complex structures on Lie algebras. 
For instance, an abelian hypercomplex structure on $\fg$, that is, is a pair of anticommuting abelian complex structures,  
 gives rise to an invariant weak HKT structure (see [\ref{cqg}] and [\ref{gp}]).

Abelian complex structures on Lie algebras were first considered in 
 [\ref{bdm}] where a construction is given  
starting with a 2-step nilpotent Lie algebra 
and applying successively a ``doubling" procedure. 
It follows from results of [\ref{ba}] that $\fa \ff \ff (\cc)$, the Lie algebra of the  affine motion group of $\cc$,  
  is the unique $4$-dimensional Lie algebra carrying 
an abelian hypercomplex structure.
  In [\ref{bar}] 
 the particular class of H-type Lie algebras was studied in detail and 
a precise answer was given to the question of  when such an algebra 
admits an abelian complex structure.

It was proved in [\ref{df}] that a real Lie algebra admitting an 
abelian complex structure is necessarily solvable. In the present article we 
give a characterization of the solvable Lie algebras admitting an abelian complex structure 
in terms of certain affine Lie algebras $\fa \ff \ff (A)$, $A$ a commutative algebra 
(Theorem~\ref{main}). These affine Lie algebras are natural generalizations  of 
$\fa\ff\ff (\cc)$ and the corresponding Lie groups are complex affine manifolds.  
 It turns out,
using the classification given in [\ref{Sn}],  that 
 all $4$-dimensional Lie algebras carrying abelian complex structures are central 
extensions of affine 
Lie algebras. 

 
In \S ~\ref{sec-obstr} we study obstructions to the existence of abelian complex structures.

\section{Complex structures on affine Lie algebras} \label{secaff}

A  complex structure on a real Lie algebra $\fg$ is an  endomorphism 
$J$ of $\fg$ satisfying 
\begin{equation}  
J^2=-Id,\quad \quad J[x,y]-[Jx,y]-[x,Jy]-J[Jx,Jy]=0, \; \;\; \forall x,y 
\in \fg.   \label{nijen}  \end{equation}
 Note that complex Lie algebras are those for which the endomorphism $J$ satisfies the stronger condition
 \begin{equation}  
J^2=-Id,\quad \quad J[x,y]=[x,Jy], \; \;\; \forall x,y 
\in \fg   \label{cx}  \end{equation}
By a hypercomplex structure we mean a 
pair of anticommuting  complex structures.

A rich family of Lie algebras carrying  complex structures is obtained by considering
  a finite dimensional real associative  algebra $A$ and $\fa \ff \ff (A)$ the Lie algebra $A \oplus A$
with Lie  bracket  given as follows:
\[   [(a,b),(a',b')]=(aa'-a'a,ab'-a'b), \hspace{1cm} a,b,a',b' \in A . \]  
Let $J$ be the endomorphism  of $\fa \ff \ff (A)$ defined by 
\begin{equation} J(a,b)=(b,-a),    \hspace{1cm} a,b  \in A .
\label{jaff1} \end{equation}
 A computation shows  that $J$ defines a complex structure 
on $\fa \ff \ff (A)$. Note that when $A$ is a vector space with the trivial product structure $ ab=0, \; a,b \in A$ one obtains   the abelian Lie algebra ${\rr}^n \oplus {\rr}^n$ with the standard complex structure $J(a,b) = (b,-a).$  Furthermore, if
 one assumes the algebra $A$ to be a complex associative  algebra,  
this extra assumption allows us to equip $\fa \ff \ff (A)$ with a pair of anti-commuting 
 complex structures.  Indeed,
 the endomorphism $K$ on  $\fa \ff \ff (A)$ 
 defined by $K(
a,b)=(-ia, ib)$ for $a,b\in A$ satisfies (\ref{nijen}) and since
$JK=-KJ$, $J$ and $K$ define a hypercomplex structure.

\begin{prop} \label{hcx} $\fa \ff\ff (A)$ carries a natural hypercomplex structure for any complex 
associative algebra $A$. 
\end{prop}

The Lie groups having Lie algebras  $\fa \ff \ff (A)$ carry invariant complex affine structures. Indeed, the bilinear map $\nabla$ given by $\nabla_{(a,b)}(c,d)= (ac,ad)$ satisfies 
\[ \nabla_{(a,b)}J(c,d)= J\nabla_{(a,b)}(c,d), \hspace{1cm} \nabla_{(a,b)}(c,d)-\nabla_{(c,d)}(a,b)= [(a,b), (c,d)] \] and $R((a,b),(c,d))= 0$ where \[ R((a,b),(c,d))= \nabla_{[(a,b), (c,d)]}-[\nabla_{(a,b)},\nabla_{(c,d)}] \]
is the curvature tensor.  In particular, using results of Boyom [\ref{Bo}],   
any such simply connected Lie group can be embedded as leaf of 
a left invariant lagrangian foliation in a symplectic Lie group.

\section{Abelian complex structures} \label{secsolv}

 An {\it abelian} complex structure on a real Lie algebra $\fg$ is an endomorphism 
of $\fg$ satisfying \begin{equation} J^2=-I, 
\hspace{1.5cm} [Jx,Jy]=[x,y], \; \;\; \forall x,y 
\in \fg.   \label{abel1}  \end{equation}
 By an abelian hypercomplex structure we mean a 
pair of anticommuting abelian complex structures. 

We  observe that  one can rewrite condition (\ref{nijen}) as follows
\begin{equation}J([x,y]-[Jx,Jy]) = [Jx,y]-[x,Jy] \; \;\; \forall x,y 
\in \fg.\label{nijen1}\end{equation} Thus, abelian complex structures are integrable. Moreover,  from (\ref{nijen1}) one has  that if 
$[x,y]-[Jx,Jy]\neq 0$ for some $x,y$    then  the commutator subalgebra has dimension $\geq 2$.  In particular, if   $\fg$ is a real Lie algebra with $1$-dimensional commutator  $[\fg ,\fg]$ then every complex structure on $\fg$ is abelian (compare with Proposition~4.1 in [\ref{bd}]).

There exist algebraic restrictions to the existence of abelian complex structures. We recall the following result
\begin{prop}[{[\ref{df}]}] \label{restr}
Let $\fg$ be a real Lie algebra admitting an abelian complex structure. 
Then $\fg$ is solvable. 
\end{prop}

Given a complex structure $J$ on a Lie algebra $\fg$,
the endomorphism $J$ extends to the
complexification ${\fg}^{\cc} = \fg \oplus i \fg$ giving a splitting
$${\fg}^{\cc} =
{\fg}^{1,0}\oplus {\fg}^{0,1} $$ where $$ {\fg}^{1,0}= \{X-iJX: X\in\fg
\} \;\; {\rm and}\;\; {\fg}^{0,1}=\{X+iJX: X\in \fg \}$$ 
are complex
Lie subalgebras of $\fg ^{\cc}$.  Using (\ref{abel1}) one verifies that abelian complex structures are those for which the subalgebras $\fg^{1,0}$ and $\fg^{0,1}$ are abelian, and conversely.  

In order to give another characterization of abelian complex structures we need first to consider  the following general class of complex structures on matrix algebras.

 Let $V$ be a real vector space, $\dim V=2n$, and fix a complex endomorphism 
 $I$ of $V$ (i.e. $I^2=-Id$).  Let us denote by $L_I$, (resp. $R_I$) 
the endomorphism of $\fg\fl(V)$ defined as $L_I(u)=I \circ u $ (resp. $R_I(u)= u\circ I$),  
$u\in \fg\fl(V)$. It is straightforward to show   that $L_I$ (resp. $R_I$) defines a 
complex structure on $\fg\fl(V)$, that is, it satisfies (\ref{nijen}). Moreover, the subalgebra $\fg\fl_{\cc}(V)$ of endomorphisms of $V$  commuting with $I$ is $L_I$ and $R_I$ invariant and the restriction of $L_I$ or $R_I$ to this subalgebra satisfies (\ref{cx}).

Consider next an arbitrary  Lie algebra $\fg$ and assume that $J$ 
is an endomorphism of $\fg$ satisfying $J^2=-Id$. In particular, 
$\dim \fg = 2n$. Consider on  $\fg^*$ the induced endomorphism, that we denote also by $J$, given by $J\alpha= -\alpha J,\; \alpha \in \fg^*$.  According to the previous observation, $R_{-J}$ is integrable on $\fg\fl(\fg)$ and $L_{J}$ is integrable on $\fg\fl(\fg^*
)$. It follows after a computation    
that $J$ is an abelian complex structure on $\fg$ if and only if the adjoint representation  $\text{ad} : 
(\fg, J)  \rightarrow (\fg\fl(\fg), R_{-J})$ is holomorphic, that is, 
$\text{ad}\,(Jx)=R_{-J} ( \text{ad}\,(x))$ for all $x\in \fg$.  Equivalently, the coadjoint representation
$\text{ad}^* : 
(\fg, J)  \rightarrow (\fg\fl(\fg^*
), L
_{J})$ is holomorphic, that is, 
$\text{ad}^*\,(Jx)=L_{J} ( \text{ad}^*\,(x))$ for all $x\in \fg$.  
 This 
paragraph can be summarized as follows: 
\begin{teo} \label{key} 
Let $J$ be a complex structure on the real Lie algebra $\fg$. Then the following conditions are equivalent:

$i)\; J$ is abelian.

$ii)$ The complex subalgebras $\fg^{1,0}$ and $\fg^{0,1}$ of  ${\fg}^{\cc}$ are abelian.  

$iii)$ The adjoint representation
ad$:(\fg , J)\rightarrow (\fg\fl(\fg),R_{-J})$ is  
 holomorphic.

$iv)$ The coadjoint representation
$\text{ad}^* : 
(\fg, J)  \rightarrow (\fg\fl(\fg^*), L_{J})$ is holomorphic.
 \end{teo}

\subsection {Examples} \label{ex3}
The simplest examples of non abelian Lie algebras carrying abelian complex structures are provided by 

\begin{enumerate}

\item [i)]$\fa \ff \ff (\rr)$, the Lie algebra of the affine motion group of $\rr$ (the bidimensional non-abelian Lie algebra), $\fa \ff \ff (\rr) = span \{ x, y\}$, with 
 bracket $[x,y]= x$  and $J$ given by $Jx=y$
and 
\item [ii)] $\rr\times \fh_n$,  where $\fh_n$ stands for the $2n+1-$dimensional Heisenberg Lie algebra,  $\rr\times \fh_n= span \{ w, z, x_i, y_i, i=1,\dots , n\}$, with non zero bracket $[x_i, y_i]= z$ and $J$ given by $Jz=w$,  
$Jx_i= y_i, i=1,\dots ,n$. 
 
\end{enumerate}

The Lie algebras introduced in i) and ii) have one dimensional commutator.  Moreover, every Lie algebra with one dimensional commutator is a trivial central extension of one of these (see Theorem~4.1 in [\ref{bd}]).   Hence we obtained the following result:  

\begin{prop} \label{1dim} Every even dimensional Lie algebra with one dimensional commutator carries an abelian complex structure. 
\end{prop}

 The next family of examples will play a crucial role in the characterization given in Theorem~\ref{main}. 

\begin{enumerate}
\item [iii)] Consider the Lie algebra $\fa \ff \ff (A)$ defined in \S \ref{secaff}
 where $A$ is a commutative algebra. 
Let 
$J$ be the complex structure on $\fa \ff \ff (A)$ defined by equation \eqref{jaff1}. 
Then one verifies  that $J$ is an abelian complex structure.  
We note that when $A={\rr}$ or $A={\cc}$, we obtain the Lie algebra of the affine motion group of  
either ${\rr}$ or ${\cc}$.
Moreover, if $A$ is a complex commutative algebra then the 
complex structure $K(a,b)=(ia,-ib)$ which anticommutes with $J$ is also abelian, hence in this case 
we obtain an abelian hypercomplex structure. 

\end{enumerate}

\begin{prop} \label{abhcx} If $A$ is a complex commutative algebra then the natural hypercomplex 
structure on $\fa \ff \ff (A)$ is abelian. 
\end{prop}

The $4$-dimensional 
Lie algebras  admitting abelian complex structures are essentially affine algebras $\fa \ff\ff (A)$ for some commutative algebra $A$ (see Proposition~\ref{dim4}). In the general situation these 
algebras are also involved as building blocks (Theorem~\ref{main}).

A particular case of the  construction just considered  occurs when one assumes  $A$ to be the   set of complex matrices of the form
\[  
  \begin{pmatrix} 0 & a_1&a_2 &\dots & a_{k-1}& a_k\\
0&0&a_1 &\dots & a_{k-2}&a_{k-1} \\
\hdotsfor[1.5]{6} \\
0&0&0&\dots & a_1&a_2 \\
0&0&0& \dots &0&a_1 \\
0&0&0 &\dots &0  &0 \end{pmatrix} .
    \]
 $A$ is commutative and 
$\fa \ff \ff (A)$ is $k$-step nilpotent, therefore existence of abelian 
complex structures imposes no restriction on the  
degree of nilpotency (compare with [\ref{df}]).

\begin{prop}
For any positive integer $k$ there exists a $k$-step nilpotent Lie 
algebra carrying an abelian hypercomplex structure.
\end{prop}

We observe that all known examples of Lie algebras carrying abelian complex structures 
are two-step solvable, but we do not know if this holds in general.

\subsection{ Main theorem}

In this section we give a characterization of solvable Lie algebras admitting abelian complex structures.  It is our aim to show that the building blocks of such algebras are the affine algebras considered in \S \ref{ex3} iii).

\begin{prop} \label{prop-affine} Let $\fs$ be a solvable Lie algebra with an abelian complex  structure $J$ admitting a decomposition $\fs=\fu + J\fu$ with $\fu$ an abelian ideal.  Then $(\fs /\fz, J)$ is holomorphically isomorphic to 
$\fa \ff \ff(A)$ for some  commutative algebra $A$.
 \end{prop}

\begin{proof} We note first that if $\fs$ is as in the statement then $\fu \cap J\fu \subset \fz$, $\fz$ the center of $\fs$.  Indeed,  if $x=Jx'\in \fu \cap J\fu$ then $[x,u]=0, u\in \fu$, and \[ [x,Ju]=[Jx',Ju]=[x',u]=0, \hspace{1cm} u\in \fu \] 
showing that $x \in \fz$.   

Let $ A= \{ ad(Jx)\; : \; x\in \fu \}$ and let $f:\fs \rightarrow \fa \ff \ff (A)$ be defined by \[ f( x+Jy)= (ad(Jy) \;, ad(Jx)).\]
  If $x'+Jy'=x+Jy$ then both, $J(x'-x)$ and $J(y'-y)$, belong to $\fz$, hence $f$ is well defined.  Clearly,  $\fz$ is contained in the kernel of $f$, since $x+Jy \in \fz$ implies that $x$ and $Jy$ are in $\fz$. Conversely, if $ad(Jy)=0= ad(Jx)$, then $x$ and $Jy$ are in $\fz$ since $J$ is abelian.   We verify next that $f$ is a Lie algebra homomorphism. If $x+Jy, \; x'+Jy'\in \fs$, then \[ f[x+Jy,x'+Jy']= (0, ad(J([x,Jy']+[Jy,x']))).\]
  On the other hand, 
\[[(ad(Jy),ad(Jx)),(ad(Jy'),ad(Jx'))]=(0, ad(Jy)ad(Jx')|_{\fu}-ad(Jy')ad(Jx)).\] 
 Now, 
\begin{eqnarray*} ad(J([x,Jy']+[Jy,x']))|_{\fu}&=& -ad([x,Jy']+[Jy,x'])J |_{\fu}\\ &=& 
ad(Jy)ad(Jx')|_{\fu}-ad(Jy')ad(Jx)|_{\fu} \end{eqnarray*}
since $ad(Jx)J|_{\fu}=0$  for $x \in \fu$, and 
\[ ad(J([x,Jy']+[Jy,x']))|_{J\fu} =0= ad(Jy)ad(Jx')|_{J\fu}-ad(Jy')ad(Jx)|_{J\fu}.
\]
Therefore, 
\[ ad(J([x,Jy']+[Jy,x']))= ad(Jy)ad(Jx')-ad(Jy')ad(Jx)
\]
showing that $f$ induces a Lie algebra isomorphism 
between $\fs /\fz$ and $\fa\ff\ff (A)$.  Moreover, $f$ is holomorphic since 
\[ fJ(x+Jy)= f(-y + Jx)= (ad(Jx) \;, -ad(Jy)).\]

\end{proof}

We show next, using a case by case verification, that the $4$-dimensional Lie algebras admitting abelian complex structures 
are fully described by the previous proposition, that is, they are central extensions of affine algebras. 
\subsection{ The 4-dimensional case} \label{sec-dim4}

The $4$-dimensional solvable Lie algebras $\fs$ carrying complex structures were classified in [\ref{Sn}] when $\dim [\fs, \fs]\leq 2$ and in [\ref{O}] when dimension of $\dim [\fs, \fs]=3$.  
From this classification one verifies that the complex structures such that $\fs^{1,0}$ and $\fs^{0,1}$ are abelian occur only when $\dim [\fs, \fs]\leq 2$ (this also follows from 
Proposition~\ref{obstr1} below). They all appear in the classification given in [\ref{Sn}] and are denoted by $S_1, S_2, S_8, S_9, S_{10},
S_{11}$ there.  We list them below:
\begin{enumerate}

\item $S_0$ : $\fs = \rr^4$.
\item $S_1$ : $\fs= \fh_1 \oplus \rr$, a direct sum of ideals, where $\fh_1$ is the 3-dimensional Heisenberg algebra (see example \ref{ex3} ii).
\item $S_2\;$: $\fs= \fa \ff\ff(\rr) \oplus \rr^2$, a direct sum of ideals.
\item $S_8\;$: $\fs= \fa \ff\ff(\rr) \oplus \fa \ff\ff(\rr)$, a direct sum of ideals.
\item $S_9\;$: $\fs= \fa \ff\ff(\rr) \oplus \rr^2$, a semidirect sum (adjoint representation) 
\item $S_{10}\;$:$\fs= \fa \ff\ff(\rr) \oplus \fa \ff\ff(\rr)$, a semidirect product of algebras (adjoint representation).
\item $S_{11}\;$: $\fs= \fa \ff\ff(\cc)$, the complexification of $\fa \ff\ff(\rr)$. 

\end{enumerate}
The above Lie algebras, modulo their center, coincide with $\fa \ff \ff (A)$ for certain commutative algebras $A$ which are listed below:
\begin{enumerate}

\item $S_0$ : $A=0$ with the trivial structure.
\item $S_1\;$ : $A= \rr$, $\rr$ with the trivial structure.
\item $S_2\; $: $A= \rr$, $\rr$ with the standard structure.
\item $S_8\;$: $A=\{  \begin{pmatrix}  
a & 0\\
0 & b 
\end{pmatrix}, a, b \in \rr \}$. 
\item $S_9\;$: $A=\{  \begin{pmatrix}  
a & 0\\
b & a 
\end{pmatrix}, a, b \in \rr \}$. 
\item $S_{10}\;$: $A=\{  \begin{pmatrix}  
a & b\\
b & a 
\end{pmatrix}, a, b \in \rr \}$. 
\item $S_{11}\;$: $A= \cc$ , $\cc$ with the standard structure.

\end{enumerate}

The above paragraph can be summarized as follows: 
\begin{prop} \label{dim4} Let $\fs$ be a $4$-dimensional Lie algebra admitting an abelian complex structure.
 Then $\fs/\fz$ is isomorphic to $\fa\ff\ff (A)$ for some commutative algebra $A$. 
\end{prop}

 The next example shows that Proposition~\ref{prop-affine} does not exhaust the class of 
Lie algebras carrying abelian complex structures. 
\subsection{Example} \label{ex} Let $\fs= \rr x_1\oplus \rr y_1 \oplus
\dots \oplus \rr x_k\oplus \rr y_k \oplus\fv$ with $\dim \fv=2n$, $k, n\geq 1$. Fix a real  endomorphism $J$ of $\fs$ such that $J^2=-I$ and 
\begin{equation}
Jx_j=y_j, \hspace{.3cm} j=1, \dots ,k,    \hspace{1.5cm} J\fv \subset \fv.
\end{equation}
Let $T_1, \dots , T_k$ be a commutative family 
of endomorphisms of $\fv$ satisfying 
\[ T_i T_j=-T_i\, J|_{\fv}\, T_j\, J|_{\fv} \hspace{.5cm} \text{for 
all }i, j .  \] 
This condition is satisfied, for instance, if $T_i$ commutes with  $J|_{\fv}$
for all $i=1, \dots , k$. 
Define a bracket on $\fs$ as follows
 \begin{equation}
[x_j,v]=T_jJv, \hspace{2cm} [y_j, v]=T_jv  \hspace{.5cm} \text{for all } v\in \fv 
\end{equation}
and extend it by skew-symmetry. 
It turns out that $\fs$ equipped with this bracket is a Lie algebra 
and $J$ becomes an abelian complex structure on $\fs$. Observe that 
$\fs$ is not in general an affine algebra, but it has the following property: there exists a $J$-stable ideal $\fs_1= \fv$ 
isomorphic to $\fa \ff\ff (\rr^n)$ such that 
$\fs /\fs_1$ is isomorphic to $\fa \ff\ff (\rr^k)$, where both, $\rr^n$ 
and $\rr^k$, are equipped with the trivial algebra structure. The general situation is described 
by the following theorem.

\begin{teo} \label{main} Let $\fs$ be a real Lie algebra and let $J$ be an abelian complex structure on $\fs$. Then there exists an increasing sequence $\{0 \}=\fs_0 \subset \fs _1 \subset 
\dots \subset \fs _{r-1}\subset \fs _r =\fs$ of $J$-stable ideals of $\fs$ such that $\fs _j /\fs_{j-1}$ is holomorphically 
isomorphic to a central extension of $\fa \ff \ff (A_j)$ with 
the abelian complex structure given by equation \eqref{jaff1},   for some commutative algebra $A_j$,  
   $1\leq j\leq r$.  

\end{teo}

\begin{proof} We proceed by induction on $\dim \fs$. 
The theorem is trivially satisfied if $\dim \fs =2$. 
If $\dim \fs >2$, we assume that the theorem is true for all Lie algebras of dimension 
strictly less than $\dim \fs$. 
Since $J$ is abelian, $\fs$ must be solvable (Proposition \ref{restr}). Let $\fu$ be a non zero abelian ideal in $\fs$, then $\fs_1= \fu +J\fu$ is a solvable Lie algebra satisfying the hypothesis of Proposition~\ref{prop-affine}.  Hence $\fs_1$ is holomorphically isomorphic to a central extension of  $\fa \ff \ff (A_1)$ with 
the abelian complex structure given by equation \eqref{jaff1}   for some commutative algebra $A_1$.  
 If $\fs_1=\fs$ we are done. Otherwise, since $\fs_1$ is a $J$-invariant ideal of $\fs$, the inductive hypothesis applies to the Lie algebra $\fs /\fs_1$ with the induced abelian complex structure. 
\end{proof}

\section{Some obstructions} \label{sec-obstr}
As a consequence of Proposition 1.5 in [\ref{ss}], if $\fn$ is a nilpotent Lie algebra 
admitting an abelian complex structure then $[\fn, \fn]$ must have codimension $\geq 3$. On the other hand, we exhibited in \S \ref{ex} solvable Lie algebras $\fs$ with 
$[\fs,\fs]$ of codimension $2k$,  $k\geq 1$, admitting abelian complex structures.    
The following result implies that if $[\fs,\fs]$ has codimension $1$ and $\dim \fs >2$ then abelian complex structures 
do not exist on $\fs$.  
\begin{prop} \label{obstr1}
Let $\fs$ be a solvable Lie algebra such that $[\fs,\fs]$ has codimension 
$1$ in $\fs$. If $\fs$ admits an abelian complex structure then $\fs$ is 
isomorphic to $\fa \ff \ff (\rr)$. 
\end{prop}
\begin{proof} Let $J$ be an abelian complex structure on $\fs$ 
and set \[
\fs= \rr a\oplus \fn
\]
where $\fn= [\fs,\fs]$ and $a$ can be chosen so that $Ja\in \fn$. Then 
\begin{equation} \begin{split}   
\fn&= \text{Im}\, \text{ad}\, (a) +[\fn ,\fn ]= 
\text{Im}\, \text{ad}\, (Ja) + [\fn ,\fn ] \\
&=\rr [a,Ja] + \text{Im}\, \text{ad}\, (Ja)|_{\fn} + [\fn ,\fn ]
=  \rr [a,Ja]  + [\fn ,\fn] \end{split} \end{equation}
and we get $\fn'\subset [\fn ,\fn']$, hence $\fn'=[\fn ,\fn']$. Now, $\fn$
is nilpotent, so we must have $\fn'=\{ 0\}$ and therefore $\fn=\rr [a,Ja]$. 
This implies the result. 

\end{proof}

As a consequence of the above proposition we obtain a large family of 
Lie algebras which do not carry abelian complex structures. In fact, consider a nilpotent Lie algebra $\fn$, $\dim \fn >1$, admitting a non-singular derivation $D$ and set 
$\fs= \rr a \oplus \fn$ where the action of $a$ on $\fn$ is given by  $D$. It follows from the proposition that there is no abelian complex structure on $\fs$. A particular case of this construction is given by Damek-Ricci extensions of H-type Lie algebras (see [\ref{dr}]). In particular,  
the solvable Lie algebras corresponding to the rank one symmetric 
spaces of non-compact type [\ref{hein}] do not admit abelian complex structures, though it is well known that they do admit complex structures (equation (\ref{nijen})).

Abelian complex structures are frequent on two-step nilpotent Lie algebras 
(see [\ref{bdm}] and [\ref{bar}]), but even in this case we have the following restriction: 
 
\begin{prop}  \label{obstr2} Let $\fn$ be a two-step nilpotent Lie algebra 
such that $2 \dim [\fn , \fn]=n (n-1)$, where $n= \dim \fn -\dim \fz\geq 3$ and $\fz$ is the center of $\fn$.  Then 
$\fn$ does not admit an abelian complex structure. 
\end{prop}
\begin{proof} We assume that $\fn$ admits an abelian complex structure $J$. 
Fix a Hermitian inner product $\langle \, ,\rangle$ on $\fn$ and consider the orthogonal decomposition $\fn=\fz \oplus \fv$. 
Being $J$ abelian, it follows that both, $\fz$ and $\fv$, are $J$-stable. 
Define a linear map $j: \fz \rightarrow$ End $(\fv)$, 
$ z \mapsto j_z$, where 
 $j_z$ is determined as follows:
\begin{equation} \label{jz}
\langle j_zv,w \rangle =
\langle z, [v,w]\, \rangle , \;\;\; \;\forall v,w \in \fv.
\end{equation} 
Observe that $j_z,\; z\in \fz$,  are skew-symmetric so that 
$z \rightarrow j_z$ defines 
 a linear map $j: \fz \rightarrow \fs \fo (\fv)$ and the restriction 
of $j$ to $[\fn , \fn]$ is injective. It follows from Lemma 1.1 
in [\ref{bar}] that $J$ commutes with $j_z$ for all $z\in \fz$, which is 
a contradiction. In fact, our assumption on $\dim  [\fn , \fn]$ says that the map  
$j: \fz \rightarrow \fs \fo (\fv)$ is surjective. Therefore, since $n\geq 3$, the only  
endomorphisms of $\fv$ commuting with all $j_z$, $z\in \fz$, are real multiples of the identity.   
\end{proof}
 Recall that a two-step nilpotent
Lie algebra $\fn$ is said to be free, of rank $n$,  when $\fz=[\fn , \fn]$ 
and $2\dim \fz =n(n-1)$, where $n=\dim \fn -\dim \fz$.   The above result says that the free two-step nilpotent Lie algebras of rank $n \geq 3$
do not admit abelian complex structures.

\section*{References}
\begin{small}
\begin{enumerate}
\item M. L. Barberis, I. G. Dotti Miatello and R. J. Miatello, 
{\it On certain locally homogeneous Clifford manifolds}, Ann. Glob. Anal. 
Geom. {\bf 13} (1995), 289--301. \label{bdm}
\item M. L. Barberis and I. Dotti Miatello, 
{\it Hypercomplex structures on a class 
of solvable Lie groups}, Quart. J. Math. Oxford 
 (2), {\bf 47} (1996), 389--404. \label{bd} 

\item M. L. Barberis, {\it Hypercomplex structures on four dimensional Lie groups}, Proc. Amer. Math. Soc. {\bf 125} (4) (1997), 1043--1054. \label{ba}
\item M. L. Barberis, {\it Abelian hypercomplex structures on central extensions of H-type Lie algebras}, J. Pure Appl. Algebra {\bf 158} (2001), 
15--23. \label{bar} 
\item N. Boyom, {\it Varietes symplectiques affines}, Manuscripta Math. {\bf 64}, (1989), 1--33. \label{Bo}
\item E. Damek, F. Ricci, {\it Harmonic analysis on solvable extensions of  H-type  groups}, J. Geom. Anal. {\bf 2} (1992), 213--248. \label{dr}
\item I. Dotti and A. Fino, {\it Hyper-K\"ahler torsion structures invariant by nilpotent 
Lie groups}, to appear in Classical Quantum Gravity. \label{cqg}
\item I. Dotti and A. Fino, {\it Hypercomplex nilpotent 
Lie groups}, preprint.  \label{df}
\item G. Grantcharov and Y. S. Poon, {\it Geometry of Hyper-K\"ahler connections with torsion}, 
Comm. Math. Phys. {\bf 213} (1) (2000), 19--37. \label{gp}
\item E. Heintze, {\it On homogeneous manifolds of negative curvature}, 
Math. Ann. {\bf 211} (1974), 23--34. \label{hein}
\item G. Ovando, {\it Invariant complex structures on  solvable real Lie groups}, Manuscripta Math. {\bf 103} (2000), 19--30. \label{O}
\item S. Salamon, {\it Complex structures on nilpotent Lie algebras}, J. Pure Appl. Algebra 
{\bf 157} (2001), 311--333. \label{ss}
\item J. E. Snow, {\it Invariant complex structures on four dimensional solvable real Lie groups}, Manuscripta Math. {\bf 66}, (1990), 397--412. \label{Sn}

\end{enumerate}
\end{small}

\end{document}